\documentclass[12pt,draftcls,onecolumn,peerreview]{IEEEtran}
\usepackage {graphicx}
\usepackage{subfigure}
\DeclareGraphicsExtensions{.eps,.pdf,.jpeg,.png}
\usepackage[cmex10]{amsmath}
\usepackage{amssymb}
\usepackage{cite}
\usepackage{bm}
\usepackage{extarrows}
\usepackage{lineno}
\title{Logic and Paradox}

\author{Xuezhi Yang, yangxuezhi@hotmail.com}

\begin{document}
\maketitle

\begin{abstract}
This article discusses the logical errors in the liar paradox, G\"odel's incompleteness theorems, Russell's paradox, and the halting problem. In order to avoid these errors, a redefinition of logic has been presented, which is concluded as four principles in set-theoretic language, including 1. Don't talk about empty; 2. Elements of set should have identity; 3. Sets should have definitions; and 4. A total set should be defined. The new definition of logic can eliminate all paradoxes and invalid statements, thus becoming a solid foundation for human language and knowledge.
\end{abstract}

\bigskip

\section{Introduction}

Logic is a set of thinking rules that all forms of life, including humans, need to abide by. Adherence to logic in thinking is essential for providing effective guidance for practices to adapt to and survive the environment.

From the birth of life, logic has been playing a fundamental role in the interaction between life and its environment. However, logic emerged as a distinct discipline not until the ancient Greek philosopher Aristotle extensively discussed what we now refer to as the three fundamental laws of logic: the law of identity, the law of non-contradiction, and the law of excluded middle  in his work "Metaphysics". In his works "Categories" and "Prior Analytics," he described the basic forms of argumentation and inference, with the syllogism as the core structure.

Logic is not only the rules of thinking in people's daily lives but also a fundamental characteristic of science. With the advancement of modern science and technology, logic has become an essential tool in stating, reasoning, and validating scientific theories, as well as in proposing, implementing, and improving various technologies.

While logic plays a role in various aspects of human life, people have also discovered paradoxical phenomena closely related to logic. Among them, the famous paradoxes include the liar paradox, G\"odel's incompleteness theorems \cite{K1931ber}, Russell's paradox\cite{Russell}, and the halting problem\cite{turing1936computable}. There are also other paradoxes, but these four are the most well-known and representative ones.

A paradox is a self-contradictory statement from which the opposite of itself can be logically derived. This raises a question: what does the existence of paradoxes indicate? Does it mean that logic has flaws? Does it mean that logic has limitations? Is there a God? All these doubts trace back to the very root of logic, mathematics, and the whole knowledge of mankind.  

If paradoxes arise even when logic is strictly followed, it does indeed indicate that there might be issues with logic.  However, if the paradoxes do not adhere to logic, the contradictions cannot be attributed to logic but rather to logical errors within the paradoxes. Upon careful examination, this article reveals that these four famous paradoxes all involve logical errors, and thus they cannot be used to demonstrate flaws in logic.

However, it is also a fact that the mainstream academic community has recognized issues with classical logic over the past century. This at least indicates that the problem of "what is logic" is still not clear. According to the current understanding, the three fundamental principles of classical logic can be summarized as follows:
\begin{itemize}
  \item Law of identity: $A = A$.
  \item Law of non-contradiction: $A\neq\neg A$.
  \item Law of excluded middle: $A$ or $\neg A$.
\end{itemize}

Among the three laws of logic, the law of identity is the most fundamental. "$A = A$"  is  explained as "Any individual entity or object is the same as itself". But this explanation is confusing and doesn't seem to convey any meaningful information. This is also why Wittgenstein asserts that self-identity is nonsense\cite{wittgenstein1921tractatus}.  As a supplementary explanation for the law of identity, it can be stated that "In the same thought process, each concept or judgment used has a definite content." Although this can be understood to some extent, the concepts used here have not been properly defined.

Logic had already been regulating the thought processes of lives even before the emergence of written language. Therefore, any textual description of logic will inevitably lead to circular definitions, making it difficult to precisely explain "what is logic." Taking this situation into account, this article carefully examines the logical errors in the four famous paradoxes and redefines the logical system, addressing the incomplete aspects of classical logic. In doing so, it avoids all paradoxes and establishes logic as a solid foundation for all human thinking and knowledge.

\section{Logical errors in paradoxes}

A paradox, by deriving its own negation from itself, violates the law of non-contradiction, which indicates the existence of  logical errors. Below, we will analyze the specific logical errors present in each of the four paradoxes separately.

\subsection{Liar Paradox}

In the 6th century BCE, philosopher Epimenides of Crete said a very famous statement: "This statement of mine is false." If Epimenides' statement is true, then its content "This statement of mine is false" must be false. On the other hand, if the statement is false, then "This statement of mine is false" is true. This leads to a contradiction, hence the name "Liar Paradox."

Epimenides' statement contains self-reference, which is expressed very ambiguously. To analyze its logical structure, it is necessary to clarify the expression.

First, let's say a regular statement, "Proposition G is false." This statement is also a proposition. If we give this proposition a name, can we call it "proposition G"? If so, proposition G becomes Epimenides'  "this statement of mine," and it asserts that itself is false, leading to self-reference and the Liar Paradox.

However, this proposition cannot be referred to as "proposition G" because it has already been used to denote some specific finished proposition before the proposition "Proposition G is false." is said.  Therefore, this proposition cannot be "proposition G" but rather a new proposition. Calling this new proposition "proposition G" would violate the law of identity.

\subsection{G\"odel's Incompleteness Theorems}

Another famous paradox is G\"odel's  incompleteness theorems. In 1900, German mathematician David Hilbert proposed a plan to establish mathematics on a fixed and complete set of axioms. In his view, all mathematical propositions should be precisely stated and proven true or false within this system. Hilbert's famous quote, "We must know, we will know," reflects his idealism in achieving this goal.
 
However, this plan was challenged by Kurt G\"odel's incompleteness theorems in 1931. G\"odel proved that in any formal system containing the axioms of natural numbers, there exist propositions that are neither provable nor disprovable within the system. He demonstrated this by constructing a "G\"odel  statement," which declares its own unprovability within the given axiom system.

The G\"odel  statement is essentially a mathematical construction of the Liar Paradox. However, if the Liar Paradox has a logical error, then G\"odel's construction process must also has a logical error. G\"odel's mistake has already been identified in \cite{yang2023}, and is reiterated in the following paragraph. 

Formula (8.1) in \cite{K1931ber} contains the expression "$subst(y,19,number(y))$". In this expression, the first parameter $y$ represents the G\"odel number of a certain proposition, and the second parameter 19 represents that the mentioned proposition contains the free variable $y$. This means that G\"odel assigns two different meanings to the symbol "$y$" and thus violating the law of identity. This error in form is exactly the same as the Liar Paradox.

Actually, "G is unprovable" is a proposition with a free variable G. If we imitate  G\"odel's error and name the proposition as "G", then we get the  G\"odel  statement in just one step and G\"odel's 20 plus pages of arguments are just totally verbose. 

\subsection{Russell's Paradox} 

Russell's paradox is a famous set-theoretic paradox discovered by the British philosopher and mathematician Bertrand Russell in 1901. Russell's paradox is regarded as  a counterexample of Cantor's naive set theory, which defines a set by the Comprehension Axiom. Cantor claimed, for any formula $\phi(x)$ containing $x$
 as a free variable, there will exist the set $\{x|\phi(x)\}$. However, Russell constructed a set as $R=\{x|x\notin x\}$, whose members are the sets that is not a member of itself. 
The question here is, whether or not $R\in R$? Just like the Liar Paradox, no matter what's the answer, it will leads to its own negation. Russell's paradox was called the third mathematical crisis, with the former two being irrational number and infinity. It is now generally accepted that  Russell's paradox is resolved by the ZFC axiomatic set theory. 

Russell's Paradox involves the question of whether a set can contain itself. We know that, except for the empty set, a set is composed of elements which need to possess identity. Identity, roughly speaking, can be described as stable determination and is the foundation of logic. We will give a more detailed and rigorous account of identity in Section III of this paper.

Forming a set requires selecting some certain elements, and then the totality of these elements becomes the set. This means the set gets its identity right after finishing selecting its elements which have had their identity already.  Therefore, the formed set must not be any of the previously selected elements. In other words, a set cannot belong to itself. 

Although a set cannot belong to itself, the proposition "for any set $x, x\notin x$" is also invalid, because when we say whether an element belongs to set $x$, that element must be an object that already possesses identity before forming set $x$.   For subject "$x$", neither "$\notin x$" nor "$\in x$" is a legitimate predicate.  

According to the currently prevailing view, Russell's Paradox shows a flaw in Cantor's naive set theory, leading to the development of the ZFC axiomatic system to address this flaw. In reality, Russell's Paradox is merely a logical error and does not indicate any flaw in naive set theory.

\subsection{Halting Problem}

Halting problem is perhaps the most well-known decision problem in computability theory. In 1936, Turing demonstrated that there can be no general procedure to decide if a computer program will eventually halt. His argument goes like this:

Suppose there exists a program Halt(P) that can determine whether any program P halts. If Halt(P) = T, it means P halts; if Halt(P) = F, it means P does not halt. Now, let's construct another program Code(P) \{ if (Halt(P) = T) loop; else return; \}. In other words, the program Code calls the program Halt, and if it detects that P halts, then Code enters an infinite loop; otherwise, Code halts. Now, what will happen if we call Code(Code)?

Through a simple examination, it can be found that if Halt detects that Code halts, then Code does not halt; conversely, if Halt detects that Code does not halt, then Code halts. This is a contradiction, by which Turing concluded that there is no program that can decide whether any program halts or not.

Similar to the liar's paradox and Russell's paradox, self-reference is involved here. As in the previous cases, these self-reference instances are logically wrong, which raises the suspicion that the Halting Problem also commits the same logical error. The crux lies in the concept of identity. 

To determine whether a program halts, the prerequisite is that the program's halting status should depend on its own rules, or has an identity, and then Halt can make the judgment. The normal call Code(P), where P is a program with identity, is logically valid. But in the self-referring call  Code(Code), Halt(Code) is called while Code is still in progress and not definite, causing a logical error. 

If Halt has the ability to detect whether or not any program will halt, it should have had the capability to detect whether a program is well defined to have identity. For an unfinished program, Halt should refuse to detect whether it will halt before  the program is finished. Thus in the self-referring call  Code(Code),  Halt(Code) will not return a T nor a F, but report an error or fall into a deadlock. If the identity detection function does not reside in Halt, the caller should ensure P is a finished program in calling Halt(P), so  Code(Code) is an illegitimate call for it fails to do that.

\section{Identity}

If one attempts to explain what logic is using words, the problem of circular definition will arise, because logic had been working with life long before language appeared. To define logic, we must start from the most basic phenomena. The words used in this section are only meant to trigger the reader's intuition. Once this intuition is triggered, these words can be discarded.

Suppose you see an apple in front of you. After looking at it for ten minutes, it remains unchanged and keep appearing as an apple. The question is: How do we explain this phenomenon? 

First, whether a thing has changed or not is a feeling, an instinct of human beings. Everyone knows whether what they sense has changed or not, this is a fundamental assumption for human being. Of course, from modern knowledge, we know that after one minute, the apple's moisture evaporates a bit, and the apple undergoes a change, but we don't feel it. This shows that human perception has a certain precision. People can only sense changes beyond a certain threshold.

Why do we need to explain this phenomenon? The purpose of explanation is to obtain an idea, which can guide human practice and ensure survival. We explain this phenomenon by stating that the apple is something with identity, something that keeps being itself, unchanged itself.  The reason we see an unchanged apple image for ten minutes is that the apple itself has not changed, and our eyes perceive it, getting a stable image of it. With this idea, we can guide our actions. We can pick it up and eat it because this idea serves as a prediction of the apple that  it will keep being an apple, something eatable and sweet, throughout the process from the moment of our decision to eat the apple to the moment we put it in our mouth. This idea enables us to survive.

Likely, if we hold an apple in our hand and release it, it will fall to the floor. Every time we do this, the same phenomenon happens. This is the identity of physical laws. Things and laws that possess identity are called objective things and objective laws. 

It cannot be ruled out that there are things in the world that do not have identity, but humans cannot comprehend these things. For example, the position of an electron, according to the uncertainty principle, does not have identity. An electron is often incorrectly described as "being both here and there", which violates the law of non-contradiction.  The correct statement is "the position of an electron doesn't have identity and can not be described". Instead, the probability of the position of an electron, which has identity,  is described in quantum mechanics in the form of the Schr\"odinger equation or the Dirac equation.

Traditionally the law of identity "$A = A$" is explained as "an entity or object  $A$ is the same as itself", which is confusing. Actually  "$A = A$" is trying to formulate identity but misses the factor of time. If we describe identity using equation, it might be "$A_{t} = A_{t+\Delta t}$, for $\Delta t < \epsilon$". It means $A$ doesn't change within $\epsilon$  so that our thinking results on $A_{t} $ are still valid for $A_{t+\Delta t}$ to guide our practices related to $A$.  The value of $ \epsilon$ depends on the aim of the practice. For example, if you want to tell somebody something about $A$, the least requirement for $ \epsilon$ is $A$ should not experience a perceivable change the moment you finish your telling.

\section{The Principle of Occam's Razor}
 
The principle of Occam's Razor was proposed by the British monk William of Ockham in the 14th century. The basic idea of Occam's Razor is "Entities should not be multiplied unnecessarily," which means that when explaining phenomena, one should prefer the simplest and most economical explanation. Therefore, it is also known as the principle of economy of thought and one of the fundamental principles of modern science. 

Just now, we used identity to explain the stable image of the apple in our sight.  Are there any other explanations? Certainly, there are. For example, the apple may not be fixed and unchanging but fluctuating in size according to some law or randomly. However, as an observer, you also fluctuate in size proportionally with the apple at the same rhythm, so the image of the apple you see appears stable. This explanation is equivalent to the one based on identity in effect, but is more complex. Therefore, according to Occam's Razor, we do not adopt this explanation.

What counts as simple or complex is also an instinct of human beings. Throughout the process of evolution, only individuals who effectively utilize resources have survived. So, the judgment of simplicity and complexity is already imprinted in the genes of life without needing an explanation. 

Because the principle of identity, which serves as the basis of logic, also relies on Occam's Razor, the principle of economy of thought is even more fundamental than logic as a guiding principle of thought.

\section{Language}

If you attempt to express and communicate what you think with others, you need language. Language serves as the tool for thought and communication, while logic provides the rules for thinking. Therefore, language needs to abide by logic.

\subsection{Goal and Verification Methods of Language}

The basic sentence structure of language is called a declarative sentence, which has the basic structure of "subject + predicate." For example, "This apple is red" is a declarative sentence, where the subject is "this apple," and the predicate is "is red." "An airplane is flying in the sky" is also a declarative sentence, with the subject being "an airplane" and the predicate being "is flying in the sky."

So, why can listeners understand what you are saying when you speak? Understanding means that after hearing what you said, the listener can reach the consensus you hoped for. The purpose of reaching consensus is for practical cooperation. Whether consensus is reached needs to be verified in practices.

For example, you say, "I have two apples here, one is red, and the other is green. Which one do you want?" The other person says, "I want the red one." After you give him the red apple, he says, "Thank you," and then starts to eat it. This is a successful cooperation that verifies you two reached a consensus.

Conversely, after you give him the red apple, he says, "I wanted the red apple, why did you give me this one?" It indicates that you two did not reach a consensus on "red", and the cooperation fails.

Once we have the goal of language and the standard for verification, we need to design the rules of language to achieve these objectives.

\subsection{Concept and Set}

From a definition perspective, a concept uses a term to refer to  a class of things; it is a set. For example, "apple," "red," and "run" are each a concept, a set. 

Proper nouns are also concepts, such as "Socrates" or "the sun," which are sets containing only one element\cite{yang2020}. 

In this context, "set" and "concept" are synonymous and can be used interchangeably. The choice of which term to use depends on word usage conventions.

\subsection{ Primitive Concept}

In language, all concepts require definitions. If we use words to define concepts, we will inevitably rely on other concepts, leading to infinite regress and circular definitions. 

In language, there is a set of earliest concepts that are not defined using words; instead, they are defined based on human instinct of categorization. These concepts are called primitive concepts\cite{yang2020}. 

For example, "apple" is a primitive concept. A kid learns this word not from a textbook or a dictionary, but from his experiences his mother gives him such a red spherical sweet and eatable thing, saying "Apple","This is an apple", "Let's have an apple, honey", "Apple time", "Would you like an apple?" and etc., in which "apple" is the common word in all these sayings. Through hundreds of times of interactions with his mother, the kid will learn the meaning of the word  "apple", because the kid has the instinct of categorization and knows an apple is quite different from a cup, a flower, a cat, a person and etc. The foundation of a person's categorization instinct  lies in the fact that the objective world is classifiable, i.e., it possesses identity. 

Other primitive concepts include  "pear" and "watermelon", "red" and "green,"  "1" and "2", and so on.  Primitive concepts form the basis of language, and all concepts defined using words ultimately rely on primitive concepts for explanation.

\subsection{Set-theoretic Expression of a Declarative Sentence}

A declarative sentence describes a property of the subject. It is worth noting that the predicate does not describe the actual situation of the subject, but rather the scope in which the actual situation of the subject falls. When we say "This apple is red," we are describing the color of this apple. But from this sentence, listeners will not know the actual color of the apple, they will only know that the actual color of the apple falls within a certain range. The word "red" represents a range among all colors. 

In the language of set theory, a declarative sentence has the form "$S\subseteq P$," where  "$S$" is the subject, and "$\subseteq P$" is the predicate,  both $S$ and $P$ are concepts or sets. We use the expression "$S\subseteq P$" rather than "$s\in P$" to represent a declarative sentence in order to cover both singular and universal propositions. 

Please note that, the sentence "This apple is red" is expressed in the language of set theory not as  \{this apple\}$\subseteq $red but as \{the color of this apple\}$\subseteq $red. Similarly, the sentence "An airplane is flying in the sky" can be expressed as \{the state of an airplane\} $\subseteq$ flying in the sky. In natural language, some omissions are made for conciseness, but these omissions should not lead to ambiguity.

\section{Principles of Logic}

\subsection{Do Not Talk about Empty}

The subject cannot be an empty set, $S\neq\phi$. Empty means there is nothing, so nothing can be said about it. If we forcefully try to speak about it, we commit a logical error. If  $S=\phi$, then both $S\subseteq P$ and $S\subseteq \neg P$ hold simultaneously, and the law of non-contradiction would not hold. 

The famous problem of empty names in the history of philosophy is about discussing nothingness. Russell argued that statements like "The present king of France is bald" and "The present king of France is not bald" are both incorrect, thus the law of excluded middle does not apply. With the rule of not talking about nothingness, "The present king of France" becomes an empty set and is not allowed to be discussed, thus resolving the problem of empty names.

\subsection{Law of Identity}

When we talk about concepts or sets, it includes the name and the content or the meaning of the concept. We use the name to distinguish different concepts, but in a context, we use the content of the concept. 

Clearly defining what a concept's name refers to and maintaining this reference relationship is called the "law of identity." In simpler terms, the law of identity means that concepts have definitions, or in other words, clear, unique, and unchanging meanings.

As we said, concept and set are synonyms. In set-theoretic rhetoric,  the law of identity means a set has definite elements.

\subsection{Chain of Identity}

The foundation of the law of identity is the identity of the objective world.  Just because concept names, i.e. linguistic symbols are also objective things which possess identity, it is possible to establish an unchanging bond between a concept name and its reference.  

Objective things have their own determinations and identity. A concept or set that takes objective things as its elements is well defined and then obtains its identity through the law of identity. This set can then serve as an element of another well defined set, which also obtains identity, so on and so forth to form a chain of identity.  All theories should be ultimately used by people to understand the world and transform the world. The chain of identity which begins with the objective world is the bridge that connects a theory and practices. 

Then we require the elements of a set have identity, either the identity of objective things, or the identity of concepts obtained from the law of identity. 

\subsection{Dimension and Total Set}

An object can be observed and described from multiple perspectives, each of which is called a dimension. For example, when describing an apple, color and weight are two different dimensions. 

In the dimension of color, if all objects had the same color, the concept of color would not have existed in human languages. The reason we have the concept of color is that objects have different colors, and we need different concepts like "red" and "green" to distinguish them.

A valid predicate on a dimension needs to be defined beforehand, and their union forms the total set of valid predicates on that dimension, denoted as $T$. 

For example, in an ancient tribe, people had only seen apples and leaves, so they used "red" and "green" to represent their colors. Then, \{red , green \} was the complete set of color predicates in this tribe. Later, they encountered pumpkins which they called them "red" under the current total set. But after a period of time, they felt that the color of pumpkins is much different from that of apples, so they decided to call the color of pumpkins "yellow." As a result, the complete set of color predicates expanded.

Before the term "yellow" was agreed upon, others would not understand the statement "The pumpkin is yellow." In other words, when we make a statement " $S\subseteq P$", there must be a pre-defined set $T$ that contains all valid predicates, satisfying $S\subset T$ and $P \subset T$. A predicate outside of $T$ is gibberish and not valid.

For a proper name $S$, which contains only one element, since it must satisfy $S\subset T$, then one and only one of  "$S\subseteq P$" and "$S\subseteq \neg P$" holds true, which corresponds to the laws of non-contradiction and excluded middle. For a group name $S$ which contains more than one elements, the negation of "$S\subseteq P$" is "There is at least one element $\alpha\in S$, $\alpha \in \neg P$". The laws of non-contradiction and excluded middle also hold. The key point here is, the definition of $ \neg P$ relies on the existence of a well defined total set $T$.

\subsection{Four Principles of Logic}

To summarize the analyses above, the logical requirements for a statement "$S\subseteq P$" are:

\begin{itemize}
  \item $S\neq\phi $;
  \item The elements of $S$ and $P$ have identity;
  \item  $S$ and $P$ have  definitions;
  \item Define a total  set $T$,  satisfying $S\subset T$ and $P \subset T$.
\end{itemize}

In comparison to traditional logical definitions, the new definition of logic has the following changes:

Since this is a textual answer to the question "What is logic," it is necessary to first explain the working principles of language from the principle of economy of thought, the world's identity, primitive concepts, etc. to the objectives, structure, and verification methods of language. This forms the foundation for answering using language, which has not been explored by previous scholars.

By clarifying that the subject cannot be an empty set and that discussions about nothingness are not allowed, the problem of empty name is resolved.

By stipulating that the elements of a set have identity, a key aspect missing in the traditional definition of logic and a critical factor leading to Russell's paradox and the halting problem, is addressed. Identity, also known as self-determination, means that an object keeps being itself and remains approximately unchanged in a short period of time. It is the most fundamental physical law, used to explain the unchanging human sensations. If the elements of a set have identity, it means they can be observed, thought about, and manipulated, leading to automatic fulfillment of the Axiom of Pairing, Axiom of Union, Axiom of Power Set, Axiom of Separation, Axiom of Replacement, Axiom of Regularity, and Axiom of Choice in the ZFC axiom system.

The law of identity is expressed as "$S$ and $P$ have definitions," meaning that these two sets have unique and clear meanings, just as the traditional interpretation of it. The law of identity serves as a rule, based on the fact that the elements of the sets and the symbols used to denote the sets have identity. The liar paradox and G\"odel's incompleteness theorems both violate the law of identity.

By defining the total set $T$ of valid predicates,  any predicate outside of $T$ is  illegitimate and the law of non-contradiction and the law of excluded middle automatically hold.

\section{Conclusions}
 
From the emergence of life to the highly developed sci-tech civilization of today, logic is omnipresent and plays a fundamental role in human thinking and practice. Logic, as the rules of thought, appeared before the advent of written language, making it a challenge to answer the question "What is logic" using words.

Aristotle's expression of logic is incomplete, which are precisely the reasons behind paradoxes. This paper discusses the liar paradox, G\"odel's incompleteness theorems, Russell's paradox, and the halting problem, and analyzes the logical errors involved. To avoid these errors, a redefinition of "what is logic" is presented.

The redefinition of logic is based on the principle of the world's identity and the principle of economy of thought. These principles rely on the instinct of life, the ability to distinguish between change and unchange, efficiency and inefficiency, which need not be expressed in words.

Language and written words are built upon the foundation of identity. The earliest concepts in language are called primitive concepts, established through human classification instinct without the need for word definitions. The concepts defined by words are ultimately based on primitive concepts.

The redefined logic consists of four fundamental principles. The first principle is that the subject cannot be an empty set, thereby avoiding the problem of empty names. The second principle is the basic principle of set theory, stating that elements must possess identity, making it the interface between all theories and the objective world. Among the traditional three logical laws, the law of identity is expressed as the third principle, while the law of non-contradiction and the law of excluded middle are replaced by the fourth principle. With the definition of a complete set of valid predicates, they become derived laws of the law of identity.

In this new logical system, all paradoxes and invalid statements are eliminated, providing human knowledge with a solid and rigorous foundation.

\bigskip\bigskip\bigskip

\bibliographystyle{IEEEtran}
\bibliography{IEEEfull,Godel}


\end{document}